\documentclass[10pt]{article}

\usepackage[applemac]{inputenc}
\usepackage{ae, aeguill}
\usepackage{graphicx}

\usepackage{amssymb}

\def \N{{\mathbb N}}

\def \R{{\mathbb R}}
\def \Z{{\mathbb Z}}

\let\io =\displaystyle
\let\± =\displaystyle

\def\A{{\cal A}}

\def\C{{\cal C}}
\def\D{{\cal D}}
\def\E{{\cal E}}
\def\F{{\cal F}}

\def\H{{\cal H}}

\def\L{{\cal L}}
\def\M{{\cal M}}

\def\V{{\cal V}}

\def\card{\mathop{\rm card}}

\def\vvar{\mathop{\rm V\!\!\!\! Var}}

\def\ccov{\mathop{\rm C\!\!\!\! Cov}}

\def\div{\mathop{\rm div}}

\def\EE{{ I\!\! E}}

\def\llbrack{\lbrack\!\lbrack}

\def\phii{\varphi}
\def\pp{{I\!\! P}}

\def\rrbrack{\rbrack\!\rbrack}

\def\tr{\mathop{\rm tr}}

\def\verti{\vert_{\infty}}
\def\Verti{\Vert_{\infty}}

\def\vvar{\mathop{\rm V\!\!\!\! Var}}

\def\notation{\stackrel {\rm def} =}

\def\bi{\bigbreak}
\def\bino{\bigbreak\noindent}
\def\sm{\smallbreak}
\def\smno{\smallbreak\noindent}

\newtheorem{thm}{Theorem}
\newtheorem{lem}{Lemma}
\newtheorem{prop}{Proposition}
\newtheorem{cor}{Corollary}

\def\squareforqed{\hbox{\rlap{$\sqcap$}$\sqcup$}}
\def\qed{\ifmmode\else\unskip\quad\fi\squareforqed}
\def\smartqed{\def\qed{\ifmmode\squareforqed\else{\unskip\nobreak\hfil
\penalty50\hskip1em\null\nobreak\hfil\squareforqed
\parfillskip=0pt\finalhyphendemerits=0\endgraf}\fi}}

\smartqed

\def\proof{\bigbreak\noindent{\it Proof.}\quad}

\catcode`\à=13 \defà{\`a}
\catcode`\á=13 \defá{\'a}
\catcode`\ä=13 \defä{\"a}
\catcode`\â=13 \defâ{\^a}
\catcode`\é=13 \defé{\'e}
\catcode`\è=13 \defè{\`e}
\catcode`\ê=13 \defê{\^e}
\catcode`\î=13 \defî{\^\i }
\catcode`\ï=13 \defï{\"\i }
\catcode`\ö=13 \defö{\"o}
\catcode`\ô=13 \defô{\^o}
\catcode`\ù=13 \defù{\`u}
\catcode`\ü=13 \defü{\"u}
\catcode`\û=13 \defû{\^u}
\catcode`\ç=13 \defç{\c c}

\begin{document}

\title{Tail estimates for homogenization theorems
 in  random media}
\author{Daniel Boivin}
\date{July 3, 2006}

\maketitle

\begin{abstract}
It is known that a random walk on $\Z^d$  among
i.i.d. uniformly elliptic random bond conductances
verifies a central limit theorem. It is also known that approximations of the covariance matrix
can be obtained by considering periodic environments.
Here we estimate the speed of convergence of this homogenization result.
We obtain similar estimates for
finite volume approximations of the effective conductance
and of  the lowest Dirichlet eigenvalue. 
A lower bound is also given for the variance of
the Green function of a random walk in a random non-negative potential.

\end{abstract}

\smno
\smno

\section{Introduction}

Consider a reversible random walk on $\Z^d$, $d\geq 1$,  with
probability transitions given by independent bond conductances, see (\ref{transpro}).
It is known from the work of Sidoravicius and Sznitman \cite{SiSz04}
that if the conductances are uniformly elliptic then
a functional central limit theorem holds. Let  $\D_0$ be the diffusion matrix.

A survey of various approximations and bounds for  $\D_0$ 
can be found in  \cite[chap. 5-7]{Hug96}
and \cite[chap. 6-7]{JKO} for this model and for related ones.
As Owhadi \cite{Owh03} showed  for  
 the jump process in a stationary random environment,
$\D_0$ can  be  approximated by the diffusion matrix  of
random walks in environments with bond
conductances  that are $N$-periodic. 
 
Bourgeat and Piatnitski \cite{BoPi04}, using results
from Yurinskij \cite{Yur86}, showed that, for a similar model,
under a mixing condition,
 the diffusion matrices converge to the homogenized matrix $\D_0$  faster than $C N^{-\alpha}$
where $C$ is a constant and $\alpha$ is a positive exponent which depends on the dimension
and the ellipticity constant.

The goal of this paper, is to obtain tail estimates for the fluctuations about the
mean of the finite volume periodic approximations and to improve 
the  estimates given in Caputo and Ioffe \cite[(1.3)]{CaIo03}.

In order to do so we apply a martingale method developed by 
Kesten in \cite{Kes93} for first passage percolation models.
This method also applies to other models where there is homogenization
and when some regularity results are available.
In all three situations that will be considered, the quantities involved are similar to first-passage
 times in that they can be expressed as solutions
of a variational problem. It is this aspect that will be exploited.

 In the second situation, tail estimates are given for
  the effective conductances of a cube. The estimates are interesting
 for dimensions $d\geq 3$. 
 Fontes and Mathieu \cite{FoMa06} considered
  random walks on $\Z^d$ with non-uniformly elliptic conductances.
  In particular, they obtained estimates on the decay of the mean return probability.
   Under similar conditions, we can prove estimates of the
 effective conductance of a cube.
   A lower bound on the variance for some distributions of the conductances
  was given by Wehr \cite{Weh97}. 
  
 In the  third situation, tail estimates are obtained
  for the spectral gap of a random walk on cubes in $\Z^d$, $d\geq 3$,
  with Dirichlet boundary conditions.

Kesten's martingale method was also used by Zerner \cite{Zer98b}
to study a random walk in a non-negative random potential.
By this method, Zerner obtained upper bounds on the variance
of the Green function. We end  this paper
with a  short calculation leading to a lower bound.

Here are some notations that will be used throughout this article.
On $\R^d$, $d\geq 1$,
the $\ell_1$-distance, the Euclidean distance and
the $\ell_\infty$-distance will respectively be denoted by
$\vert\cdot\vert_1$, $\vert\cdot\vert$ and $\vert\cdot\vert_\infty$.
$\eta\in\R^d$ will be considered as a column vector and its transpose will
be denoted by $\eta'$ so that $\eta^2=\tr (\eta\eta' )$.
We say that two vertices of $x,y\in\Z^d$ are neighbours, and we will write 
$x\sim y$, if $\vert  x-y\vert=1$.
For $u:\Z^d\to\R^d$, let $\Vert u \Verti = 
\sup_{x\in\Z^d}\vert u(x)\verti$.
For $m < n \in\N$, $\llbrack m,n\rrbrack = \lbrack m,n\rbrack\cap \N$.
With this notation, for an integer $N\geq 1$, $Q_N$ will be the cube in $\Z^d$
defined by
$Q_N=\llbrack 1,N\rrbrack^d$,   $\overline Q_{N}=\llbrack 0,N+1\rrbrack^d$
and  its boundary is $\partial Q_N= \overline Q_{N}\setminus Q_{N}$.

\section{The stationary environment}

Let $ a(x,y,\omega)$, $x,y\in\Z^d$, $x\sim y$ be a sequence
of random variables on a probability space
$(\Omega,\F,\pp)$ such that  $a(x,y,\omega)=a(y,x,\omega)$ for all  $x\sim y$.
The random variable $a(x,y,\omega)$ can be interpreted as the electric or thermic
conductance of the edge joining $x$ and $y$.

We will assume that this sequence is stationary. That is,
there is a group of measure preserving
transformations $(T_x ; x\in\Z^d)$ acting on $(\Omega,\F,\pp)$
such that for all $x\sim y$ and  $z\in\Z^d$, $a(x+z,y+z,\omega)=a(x,y,T_z\omega)$.
The expectation with respect to $\pp$ will be denoted by $\EE$ or by $\langle \cdot\rangle$.

Let $\±\L^d$ be the set of edges in $\Z^d$.
In an environment $\omega$, the conductance of an edge $e\in\L^d$
with endpoints $x\sim y$ will be denoted by $a(x,y,\omega)$ or by $a(e,\omega)$.

For most results, we will also assume that 
the conductances are uniformly elliptic: 
there is a constant $\kappa\geq 1$, called the ellipticity constant, such that for all
$x,y\in\Z^d$, $x\sim y$,
and $\pp$-a.s.,
$$\kappa^{-1}\leq a(x,y,\omega)\leq \kappa.$$

Given an environment $\omega\in\Omega$,
let $(X_n;n\geq 0)$ be the reversible random walk on
$\Z^d$ with transition probabilities given by
\begin{equation}\label{transpro}
p(x,y,\omega) = a(x,y,\omega)/a(x,\omega),\quad x\sim y,
\end{equation}
where $\io a(x,\omega) =\sum_{y\sim x} a(x,y,\omega)$
is a stationary measure.
These transition probabilities induce a probability $P_{z,\omega}$ on the paths of the
random walk  starting at $z\in\Z^d$. The corresponding expectation will be denoted by
$E_{z,\omega}$.

\bi
The following proposition can be found under various forms in 
\cite{BoSz02}, \cite{JKO},  \cite{KiVa86}, \cite{Koz85} among  others. 
It can be shown using
Lax-Milgram lemma and Weyl's decomposition.
The corrector field can also be constructed directly
using the resolvent of the semigroup.
 
 \begin{prop}
Let $(a(e) ; e\in\L^d)$ be a stationary sequence of uniformly elliptic conductances.
Then there is a unique function 
$\± \chi_0 :\Z^d\times\Omega \to \R^d,$
called the corrector field of the random walk, which verifies
\begin{enumerate}
\item\label{integra}
$\±\chi_0(x,\cdot)\in(L^2(\pp))^d$ for all $x\in\Z^d$,

\item
$\±\int_{\Omega} \chi_0(x,\omega)\pp(d\omega) =0$ for all $x\in\Z^d$,

\item\label{cocycle}
 the cocycle property : for all $x,y\in\Z^d$ and $\pp$-a.s.
$$\chi_0(x+y,\omega) = \chi_0(x,\omega) +\chi_0(y,T_x\omega).$$

\item\label{harmocoord}
the Poisson equation : 
for all $x\in\Z^d$ and $\pp$-a.s.
$$E_{x,\omega}(X_1 + \chi_0(X_1)) = x+\chi_0(x).$$
\end{enumerate}

\end{prop}

The last property shows that, $\pp$-a.s.,
$\± X_n + \chi_0(X_n)$, $n\geq 0$,
is a martingale under $P_{0,\omega}$.
This martingale and the corrector field are carefully investigated
 in \cite{SiSz04}  to prove that, $\pp$-a.s., 
the reversible random walk starting at the origin verifies a functional
central limit theorem. In particular, it 
verifies a central limit theorem with a covariance matrix given by
\begin{equation}\label{statcov}
\D_0   =  \langle a(0)\rangle^{-1}\int\sum_\Lambda a(0,z)(z+\chi_0(z))(z+\chi_0(z))'d\pp.
\end{equation}
Note that  in a stationary environment $\D_0$ might not be a diagonal matrix.

\section{The periodic approximation}

Given an environment $\omega\in\Omega$
and an integer $N\geq 1$, introduce an environment $N$-periodic  on $\Z^d$,
$d\geq 1$,  by setting
$$\dot a_N(x,y,\omega)=a(\dot x,\dot y,\omega),\quad x\sim y$$
where $\dot x, \dot y \in \llbrack 0,N\rrbrack^d$, $\dot x\sim\dot y$ and
$\dot x\equiv x$, $\dot y\equiv y$ mod $N$ coordinatewise.

Then consider the reversible random walk on $\Z^d$
with transition probabilities given by
  $$\dot p_{N}(x,y,\omega)=\dot a_N(x,y,\omega)/ \dot a_N(x,\omega),\quad x\sim y$$
   where $\dot a_N(x,\omega)=\sum_{y\sim x}\dot a_{N}(x,y,\omega).$
 These  induce a probability $\dot P_{z,N,\omega}$ on the paths  starting at $z\in\Z^d$.
  The corresponding expectation will be denoted by $\dot E_{z, N,\omega}$.
  We will also use the Laplacian $\dot H_{N,\omega}$ which
  is defined on the set of functions $u : \Z^d\to\R$ by
  $\dot H_{N,\omega}u(x)= u(x) - \dot E_{x, N,\omega}u(X_1)$.
  
\subsection{Periodic corrector fields}\label{regres1}

As it was done for stationary environments, it is important 
to construct a periodic corrector field $\dot\chi_N : \Z^d\times \Omega \to\R^d$ 
for the random walk $(X_n;n\geq 0)$ in the periodic environment
so that $\pp$-a.s.
$$X_n +\dot\chi_N(X_n),\qquad n\geq 0,$$
is a martingale with respect to $\dot P_{0,N,\omega}$.

Therefore, $\pp$-a.s., $\dot\chi_N$ must verify  the equations 
$\± \dot E_{x,N} (X_1 +\dot\chi_N(X_1)) = x +\dot\chi_N(x)$
for all $x\in\Z^d$, or equivalently,
\begin{equation}\label{drift}
   \dot H_{N} \dot \chi_N(x) = \dot d_N(x),\quad \hbox{for all} \   x\in\Z^d,
\end{equation}
where $\± \dot d_N(x) =\dot E_{x,N}(X_1) -x$ is the drift of the walk.
Note that  each coordinate of $\dot d_N$
is  $N$-periodic.

The  vector space of $N$-periodic functions on $\Z^d$ can be identified with
$$\dot\H_N =\{ u:\overline Q_N \to \R\   ;\   u(x)=u(y)\   \forall x,y\in\overline Q_N
\hbox{  such that   }\   x\equiv y\  {\rm   mod }\   N\}.$$
Note that if $x\equiv y$ mod $N$ then $x$ or $y\in \partial Q_N$.
Each function $v: Q_{N}\to\R$ 
has a unique extension to $\partial Q_{N}$ which belongs to
$\dot\H_N$. For $u\in\dot\H_N$, define $\dot H_{N} u$ on
$\partial Q_{N}$ so that it belongs to $\dot\H_N$. Then
 $\dot H_{N}: \dot\H_N\to\dot\H_N$ is  a bounded linear operator.
 
  For two functions $u,v:\overline Q_N\to\R$, define
  the norm, the scalar product,  and the Dirichlet form respectively by
  $\±\Vert u\Vert^p_{p,\dot N} = \sum_{x\in Q_N} \vert u(x)\vert^p \dot a_N(x)$, $1\leq p<\infty$, 
 $\± (u,v)_{\dot N} =\sum_{x\in Q_N} u(x)v(x) \dot a_N(x)$ and
   \begin{eqnarray*}
      \dot\E_N(u,v) &  =  & 
\sum_{x,y} \dot a_N(x,y)  (u(x)-u(y))(v(x)-v(y))
  \end{eqnarray*}
 where  the sum is over all ordered pairs $\{x,y\}$
 such that $x\in Q_{N}$ and $y\in\llbrack 0, N\rrbrack^d$.

This expression makes sense for all functions $u,v :\Z^d \to \R$.
But if both $u,v$ are $N$-periodic, then the Green-Gauss formula holds
\begin{equation}\label{GG}
    \dot\E_N (u,v) = (u, \dot H_N v)_{\dot N}.
\end{equation}

\sm
Let  $\±\dot\H_N^0(\omega)=\{ u\in\dot\H_N\  ;\   (u, 1)_{\dot N} = 0\}.$
Note that $\io \dot\H_{N}^0$ depends on the environment but 
$\io\dot\H_{N}$  does not.

All the properties of the solutions of a Poisson equation that will be needed
are gathered in the following proposition.

\bi
\begin{prop}\label{poisson}
Let $(a(e) ; e\in\L^d)$, $d\geq 1$,  be a stationary sequence of uniformly elliptic conductances.
Then $\pp$-a.s., 
\begin{enumerate}
\item\label{hnivert}
$\± \dot H_{N} : \dot\H_N^0\to \dot\H_N^0$ is a  bounded invertible linear operator .

\item\label{poissonex}
For $f\in\H_N$,
$\dot H_{N} u = f$ possesses a solution $u\in\dot \H_{N} $
if and only if $f\in \dot\H_{N}^0$. 

\item 
Let  $f\in \dot\H_{N}^0$.
\begin{enumerate}
\item \label{varpb1}
The infimum 
$\± \inf_{\dot\H_{N}^0} \lbrack \dot\E_N(u,u) -  2(u,f)\rbrack $
is attained by the solution 
 $u\in\dot \H_{N}^0 $ of the equation
$\dot H_{N} u =  f$.

\item\label{varpb2}
The infimum  $\± \gamma := \inf_\M \dot\E_N(u,u)$
where $\± \M =\{ u\in\dot\H^0_N\  ; \  (f,u)_{\dot N}=1 \}$
is attained by the solution 
 $u\in\dot \H_{N}^0 $ of the equation
$\dot H_{N} u = \gamma f$.
\end{enumerate}

\item\label{pardoux}
If $f\in \dot\H_{N}^0$ then the unique solution $u\in\dot\H_N^0 $ of 
$\dot H_{N} u = f$ is
$$u=\int_0^\infty e^{-t\dot H_N}f dt,\quad  x\in Q_N.$$

\item\label{regularity}
There is a  constant  $C=C(d,\kappa)<\infty$ such that
for  all $N\geq 1$ and $f\in \dot\H_{N}^0$,  $u\in\dot\H_N^0$,  the  solution of $\dot H_{N} u = f$,
verifies the regularity estimates,
$$\Vert u\Vert_\infty \leq C N^2 \Vert f\Vert_{2,\dot N}\quad 
\hbox{ and  }\quad
\Vert u\Vert_\infty \leq C N^{2}(\log N) \Vert f\Vert_\infty.$$
\end{enumerate}
\end{prop}

\proof
For all $u\in\dot\H_{N}^0$, $\dot H_{N}u\in\dot\H_{N}^0$  
by the Green-Gauss formula (\ref{GG}).

$\dot H_{N}:\dot\H_{N}^0\to\dot\H_{N}^0$ is invertible since if 
$\dot H_{N}u=0$ then $u$ is constant by the maximum principle.

The variational principle \ref{varpb2}  holds for the Poisson equation on a smooth compact
Riemannian manifold with $f\in\C^\infty$, see
\cite[proposition 2.6 due to Druet]{Heb00}. The same arguments can be used.

Suppose that $f$ is not identically 0.
Then $\M$ is a closed convex set which is not empty since $\Vert f\Vert^{-2}_{2,N} f \in \M$.
Therefore the infimum,  $\gamma$,  is attained for some $u_0\in\M$
and $\gamma >0$ since $u_0$ is not constant.
Then using a theorem by Lagrange,
 there are two constants $\alpha,\beta\in\R$ such that
for all  $x\in Q_N$,
$$2\sum_{y\sim x}(u_0(x) - u_0(y))\dot a_N(x,y)
 - \alpha f(x) \dot a_N(x) -\beta \dot a_N(x) = 0$$
and  $\± 2\dot H_N u_0(x)  - \alpha f(x)  -\beta  = 0.$
 Therefore, for all $\± \phii\in \dot\H_N$,
 $$ 2\dot \E_N(u_0 ,\phii) 
 =\alpha(f,\phii)_{\dot N} + \beta (1,\phii)_{\dot N}.$$
 For $\phii=1$, one finds that $\beta =0$ while for $\phii=u_0$,
 one finds that $\alpha=2\gamma$.
Hence $(Hu_0,\phii)_{\dot N}  = \gamma(f,\phii)_{\dot N} $ for all 
$\± \phii\in \dot\H_N$.
  
The variational principle \ref{varpb1} can be proven similarly.

To prove the last two properties, the estimate of the speed of convergence
to equilibrium of a Markov chain on a finite state space given in terms of the spectral gap
is needed.
See for instance  \cite[Section 2.1]{Sal97}.

Let $\dot K_{N}(t,x,y)$, $t\geq 0$ and $x,y\in Q_N$,
be the heat kernel of $\io e^{-t\dot H_{N}}$.
Then for all $t\geq 0$ and $x,y\in Q_N$,
$\dot K_{N}(t,x,y) \geq 0$ and $\± \sum_y \dot K_{N}(t,x,y) = 1$.
In particular, for all $f\in\dot\H_N$ and $t\geq 0$,
\begin{equation}\label{infcontract}
\Vert e^{-t\dot H_{N}} f\Verti\leq \Vert f\Verti.
\end{equation}

Denote the volume of the torus $Q_N$,
 the invariant probability 
for the random walk on  $ Q_N$ and 
the smallest non zero eigenvalue
of $\dot H_N$ on $\dot\H_N$ respectively by
\begin{eqnarray*}
\dot a_N(Q_N)=\sum_{x\in Q_N} \dot a_N(x),\quad
\dot\pi_{N}(x)= \dot a_N(x)/\dot a_N(Q_N)
\quad\hbox{and}\quad \dot \lambda_{N} .
\end{eqnarray*}

Then for all $t>0$ and $x,y\in Q_N$,
\begin{equation}\label{expdecay}
 \vert \dot K_{N}(t,x,y) - \dot \pi_N(y)\vert
\leq \kappa \exp\big( -t\dot \lambda_{N}\big).
\end{equation}

Therefore, for all $t>0$ and  $f\in \dot H_{N}^0$,
\begin{equation}\label{1contract}
\Vert e^{-t\dot H_{N}}f\Verti = \sup_{x}\vert \sum_y (\dot K_N(t,x,y) - \dot\pi_N(y)) f(y)\vert
\leq \kappa \exp(-t\dot \lambda_{N}) \Vert f\Vert_{1,\dot N}
\end{equation}
By the Riesz-Thorin interpolation theorem, from (\ref{infcontract})  and  (\ref{1contract}),
we obtain that 
\begin{equation}\label{2contract}
\Vert e^{-t\dot H_{N}}f\Verti \leq \sqrt\kappa 
\exp(-t\dot \lambda_{N}/2) \Vert f\Vert_{2,\dot N}
\end{equation}

This will be completed by the following lower bound on $\±\dot \lambda_{N}$.
There exists a  constant $C_{1}>0$ which depends only on the dimension and on the
ellipticity constant $\kappa$ such that $\pp$ a.s. and
for all $N\geq 1$,
\begin{equation}\label{lowerbdest}
N^2\dot \lambda_{N} > C_{1}.
\end{equation}

This follows from the Courant-Fischer min-max principle 
\cite[p. 319]{Sal97} by comparison with the eigenvalues
of the simple symmetric
 random walk which corresponds to the case where the conductance
of every edge is $1$. For Neumann boundary conditions
the expressions are not as explicit but for Dirichlet and periodic 
boundary conditions on $Q_{N}$, the eigenvalues 
can be calculated explicitely much as in \cite{Spi64}. 
We find that for each $\xi\in \llbrack 0,N\llbrack^d$, 
there is an eigenvalue for the periodic boundary conditions on $Q_{N}$, 
$\lambda_{\xi}(Q_{N})$, that verifies
$$\lim_{N\to\infty} N^2\lambda_{\xi}(Q_{N}) = \frac{\pi^2}{d }\sum_{\vert 
z\vert=1}(\xi\cdot z)^2\quad\hbox{as}\quad N\to\infty.$$

The representation formula given in \it \ref{pardoux}\rm\    follows from the spectral estimates
(\ref{expdecay}) and (\ref{lowerbdest}). See \cite{PaVe01} for another recent application of 
\it \ref{pardoux}.\rm

The first regularity result follows from the representation formula (\ref{pardoux})
and (\ref{2contract}) : for $f\in\dot\H_N^0$,
\begin{eqnarray*}
\Vert u\Verti
  \leq   \int_0^\infty \Vert e^{-s\dot H_{N}}f\Verti ds 
  \leq   \int_0^\infty  \sqrt \kappa  e^{-s \dot \lambda_{N}  /2} \Vert f\Vert_{2,\dot N} ds 
  \leq C N^2 \Vert f\Vert_{2,\dot N}.
 \end{eqnarray*}
 The second one follows from (\ref{infcontract}) and (\ref{1contract}) :
for $f\in\dot\H_N^0$ and $t>0$, 
\begin{eqnarray*}
\Vert u\Verti
  \leq  \int_0^t \Vert e^{-s\dot H_{N}}f\Verti ds + \int_t^\infty \Vert e^{-s\dot H_{N}}f\Verti ds
  \leq  t \Vert f\Verti +  \frac{ e^{- t\dot \lambda_{N}}}{\dot \lambda_{N}} \Vert f\Vert_{1,\dot N}
 \end{eqnarray*}
Use $ \Vert f\Vert_{1,\dot N} \leq 2d\kappa N^d \Vert f\Verti$
and let $\± t=\frac{d}{C_1}N^{2}\log N$.
\qed
\bi
For a function $u:\overline Q_N\to\R^d$,
define $\dot H_N u$ in $Q_N$ by applying it coordinatewise.
 \sm
Let $g:\Z^d\to \Z^d$ be the function defined by $g(x) = x$.
\bi
\begin{cor}\label{defchin}
Let $(a(e) ; e\in\L^d)$, $d\geq 1$,  be a stationary sequence of uniformly elliptic conductances.
Then for  all $N\geq 1$ ,
there is a unique function $\±\dot \chi_N : \overline Q_N \times \Omega \to \R^d$
such that   $\pp$-a.s.,
in each coordinate, it is  in $\dot \H_N^0(\omega)$ and
 \begin{equation}\label{poisseq}
\dot H_{N}\dot\chi_N =  -\dot H_{N} g,\quad\hbox{on}\  Q_N.
\end{equation}

Moreover, there is a  constant $C = C(d,\kappa)$ such that
 \begin{equation}\label{bdchin}
\Vert \dot\chi_N\Vert_\infty \leq C N^{2}\log N.
\end{equation}
\end{cor}

\proof
Note that for each coordinate of 
$ -\dot H_{N} g$ belongs to  $\dot\H_N^0$ :

Indeed  $\dot H_{N} g$ is $N$-periodic and 
\begin{eqnarray*}
\sum_{x\in Q_N}\dot  H_{N}g(x)\dot a_N(x) 
 & = &  \sum_x\sum_{y\sim x}\dot a_N(x)\dot p_N(x,y)(g(x)-g(y))\\
  & = &  \sum_x\sum_{y\sim x}\dot a_N(x,y) (x-y) = 0.
\end{eqnarray*}

\sm
Then use   proposition \ref{poisson} for the function $\± f=  -\dot H_{N} g$.
The regularity estimate (\ref{bdchin}) follows from property \ref{regularity}
since $\Vert f\Verti \leq 1$.
\qed

\bi
The next step is to express the covariance matrix of the walk 
in a periodic environment in terms of $\dot\chi_N$.
By (\ref{drift}),  $\io M_n =X_n + \dot\chi_N(X_n)$, $n\geq 0$ is a martingale
with uniformly bounded increments :  $\±  Z_n=M_n - M_{n-1}.$

Let  $\± h(x) = \dot E_{x, N} (Z_1Z'_1)$.
Since $\±  h\in\dot\H_N$,
  by the ergodic theorem for a Markov chain on $Q_N$, $\dot P_{0,N}$ a.s.,
\begin{eqnarray*}
   \frac{1}{n}\sum_1^n \dot E_{0,N}(Z_j Z'_j\mid X_{j-1})
  &  = &  \frac{1}{n}\sum_1^n h( X_{j-1})
   \to  \sum_{Q_N} \dot \pi_N (x) h(x)\   \hbox{as}\   n\to \infty. \\
  \end{eqnarray*}

Then by the martingale central limit theorem
(see \cite[(7.4) chap. 7]{Dur91}),
$\±\frac{1}{\sqrt n}M_n$
converges to a Gaussian law. Hence  $\±\frac{1}{\sqrt n}X_n$ also converges
to a Gaussian with the same covariance matrix which is given by
\begin{equation}\label{dnchin}
 \dot \D_N =  \sum_{Q_N} \dot \pi_N (x) h(x)
\end{equation}
\begin{eqnarray*}
= \dot a_N(Q_N)^{-1} \sum_{x\sim y}
\dot a_N(x,y)(\dot v_N (y)-\dot v_N(x)) ( \dot v_N(y)-\dot v_N(x))'.
\end{eqnarray*}
where $\± \dot v_N(x) = x  + \dot\chi_N(x)$.

For uniformly elliptic, stationary and ergodic conductances,
Owhadi \cite[theorem 4.1]{Owh03} showed that 
for the jump process
$\dot\D_{N}$, the effective diffusion matrix in the periodic
environment, converges to the homogenized effective diffusion matrix $\D_{0}$ .
And since the jump process and the random walk on $\Z^d$
have the same diffusion matrix, the convergence theorem holds:
$\pp$-a.s. as $N\to\infty$,
\begin{equation}\label{convdifma}
\dot\D_N\to\D_0.
\end{equation}

 It is shown in \cite{Owh03}  using the continuity of
Weyl's decomposition and in \cite{BoPi04} for a diffusion with random 
coefficients. They are both illustrations of the principle of periodic 
localization \cite[p. 155]{JKO}.
At  the end of this section, a terse proof
by homogenization is given.

\bi
To write $\dot\D_N$ in terms of  the Dirichlet form on $\dot\H_N$,
we will extend  the definition of $\dot\E_{N}$ to $\R^d$-valued functions 
so that the expression of $\dot\D_{N}$  given in (\ref{dnchin})
becomes
\begin{eqnarray*}
\dot\D_{N} = \dot a_{N}(Q_{N})^{-1}\dot\E_{N}(\dot v_{N}, \dot v_{N})
\end{eqnarray*}
where $\± \dot v_N =   g +\dot\chi_N $.

For two functions $u,v:\overline Q_N\to\R^d$, 
 define
  \begin{eqnarray*}
  (u,v)_{\dot N} & = &  \sum_{x\in Q_N} u(x) v(x)' \dot a_N(x),  \\
    \hbox{and}\quad  \dot\E_N(u,v) &  =  & 
\sum_{x,y}\dot a_N(x,y)  (u(x)-u(y)) (v(x)-v(y))'
  \end{eqnarray*}
 where  the sum is over all ordered pairs $\{x,y\}$
 such that $x\in Q_{N}$ and $y\in\llbrack 0, N\rrbrack^d$.
 
\bi
Coordinatewise,
the periodic corrector fields are the solutions of
variational problems.
Indeed, from the variational formula in \ref{varpb1} of proposition \ref{poisson}, 
we have that
$\pp$ a.-s.  and  for all $N\geq 1$,

\begin{eqnarray*}
 \inf\{\tr\dot\E_N(g+u&,&  g+u) ; u \in  (\dot\H_N^0)^d\}\\
 & = &  \tr\dot\E_N(g , g) + \inf\{\tr\dot\E_N(u , u)-2 \tr (f , u)  ; u \in  (\dot\H_N^0)^d\}\\
   & =  & \tr\dot\E_N(\dot v_N,\dot v_N)
\end{eqnarray*}
where $g(x) = x$ and $f =-\dot H_N g$ as in corollary \ref{defchin}.
  
In particular, since
$$  \tr  \dot\E_{N}(\dot\chi_{N},\dot\chi_{N} )
\leq 2   \tr  \dot\E_{N}(\dot v_{N},\dot v_{N} )
+2  \tr  \dot\E_{N}(g,g ) \leq 4  \tr  \dot\E_{N}(g,g ),
$$
there is a constant $C=C(d,\kappa)<\infty$ such that
 $\pp$-a.s. and for all $N\geq 1$,
\begin{equation}\label{chinln}
   \tr  \dot\E_{N}(\dot\chi_{N},\dot\chi_{N} )\leq C N^d.
   \end{equation}
   
   The second variational principle, \ref{varpb2} of proposition \ref{poisson},   
   could be used to obtain a lower bound on 
   $\±  \tr\dot\E_N(\chi_N , \chi_N) $.
         
\subsection{Further regularity results}\label{regres2}

In the following proposition, we improve  the estimate given in 
  (\ref{bdchin}) for dimensions  $2\leq d\leq 4$.

\begin{prop}
Let $(a(e) ; e\in\L^d)$, $d\geq 2$,  be a stationary sequence of uniformly elliptic conductances.
 Then there is a constant $C=C(d,\kappa)<\infty$,
such that for all $N\geq 1$
\begin{equation}\label{bdchin24} 
\Vert \dot\chi_{N}\Verti\leq \left\{\matrix{
    C N^{d/2}  & \hbox{for $d\geq 3$} \cr
   CN(\log N)^{1/2} & \hbox{for $d = 2$} \cr}
\right.
\end{equation}
\end{prop}

\proof
For $\eta\in\R^d$, let $z_0 $ and $z_1$ be two
vertices of $\Z^d$ such that
$$\eta\cdot\dot\chi_{N}(z_{0}) = \min_{x\in\Z^d}
\eta\cdot\dot\chi_{N}(x)
\quad\hbox{and}\quad 
\eta\cdot\dot\chi_{N}(z_1 ) = \max_{x\in\Z^d} 
\eta\cdot\dot\chi_{N}(x).$$
\sm
Since $\dot\chi_{N}$ is $N$-periodic, we can assume that
$\vert z_0 - z_1\verti \geq N$.
\sm
Let $z_0 =x_{0},x_{1},\ldots,x_{n-1},x_{n}= z_1 $
be a path from $z_{0}$ to $z_1$
such that $1\leq n\leq dN$.
\sm
For $i=0$ or $1$, let
$\io \overline w_{i}(y) =\left( \frac{\vert z_1- z_{0}\verti}
{1+\vert y -z_{i}\verti} \right)^{d-1}$ and
$$ {\cal P}_{i} =\{ y \in \Z^d ;\vert y -z_{i}\verti \leq 
\frac{1}{2} \vert z_{1} -z_{0}\verti\}.$$
\sm 
For $\C$, a given set of finite paths in $\Z^d$,
let $\io w(y)=\card\{ \gamma\in\C ; y\in\gamma\}$.
\sm
Then by \cite[lemma 2, p.26]{BoDe91},
there is a constant $C<\infty$, which depends only on the 
dimension $d$,
and there is a set  of paths $\C$ from $z_{0}$ to $z_1 $ such that
$$\card \C = \vert  z_1 - z_0 \verti^{d-1}$$
and  such that for all $y\in\Z^d$,
\begin{equation}\label{poids}
w(y) \leq  C\   \overline w_{i}(y)  \   \hbox{ if $y\in{\cal P}_{i} $} \  \hbox{and   $w(y)=0$  otherwise.}
\end{equation}

For each path of $\C$, $z_0 =x_{0},x_{1}\ldots,x_{n-1},x_{n}= z_1 $
and for all $x\in\Z^d$,
$$\vert \eta\cdot\dot\chi_{N}(x)\verti 
    \leq  \vert \eta\cdot \dot\chi_{N}(z_1) -  \eta\cdot \dot\chi_{N}(z_{0})\verti
\leq\vert \eta\verti \sum_{j=1}^n \vert  \dot\chi_{N}(x_{j})- \dot\chi_{N}(x_{j-1})\verti 
$$
\sm
Since this holds for all paths in $\C$, it also holds for the 
arithmetic average over the paths of $\C$. Therefore,
$$\vert \eta\cdot \dot\chi_{N}(x)\verti 
\leq \frac{\vert \eta\verti}{\vert z_1 - z_0 \verti^{d-1}}\sum_{y\in {\cal P}_{0}\cup{\cal P}_{1}} w(y)  h_{N}(y)
$$
where $\io h_{N}(y) = \sum_{z\sim y }\vert \dot\chi_{N}(y) - 
\dot\chi_{N}(z)\verti$.
\sm
By  (\ref{chinln}) and  (\ref{poids}), 
$\±  \vert z_1 - z_{0}\verti^{1-d}\sum_{y\in 
{\cal P}_{1}} w(y)  h_{N}(y) $
\begin{eqnarray*}
& \leq & C \vert z_1 - z_{0}\verti^{1-d}\sum_{y\in 
{\cal P}_{1}} \overline w_{1}(y)  h_{N}(y)\\
& \leq & \frac{C}{N^{d-1}}\left(\sum_{y\in 
{\cal P}_{1}} \overline w_{1}^2(y)\right)^{1/2}
\left(\sum_{y\in 
{\cal P}_{1}} h_{N}^2(y)\right)^{1/2}\\
& \leq & 
\frac{C}{N^{d-1}}\left(\int_{1}^{N}\left(\frac{N}{r}\right)^{2(d-1)}r^{d-1}dr\right)^{1/2}
\left(\tr  \dot\E_{N}(\dot\chi_{N},\dot\chi_{N})\right)^{1/2}\\
 & \leq & 
C N^{d/2}  \left(\int_{1}^{N}r^{1-d}dr\right)^{1/2}
\end{eqnarray*}
where the constant $C$ now depends  on $\kappa$ and $d$.

And similarly for the sum over ${\cal P}_{0}$.
\qed
 
\bi
In the next section the $L^\infty$-estimates  (\ref{bdchin}) 
and (\ref{bdchin24}), will be combined with the following Hölder regularity result
shown in  \cite[prop. 6.2] {Del97} by J. Moser's iteration method
for reversible random walks on  
infinite connected locally  finite graphs
 with uniformly elliptic conductances :
\smno\it
Let $(a(e) ; e\in\L^d)$, $d\geq 1$,  be a (non-random) sequence of uniformly elliptic conductances.
Then there are constants $\alpha>0$ and $C<\infty$, 
which depend only on the dimension and on the ellipticity constant,
such that if for  $N\geq 1$, 
 $\io u:\overline Q_{2N} \to\R$ verifies  $\dot H_{2N} u=0$  in $Q_{2N}$,
then for all $x,y\in Q_{N}$,
\begin{equation}\label{holdreg}
    \vert u(x) - u(y)\vert 
< C\left(\frac {\vert x-y\verti }{  N}\right)^\alpha\max_{Q_{2N}} \vert 
u\verti.
\end{equation}
\rm\  

\subsection{A proof of (\ref{convdifma}) by homogenization}

This is analogous to the
 homogenization results of \cite[chapters 7 - 9]{JKO}.
An appropriate framework for random walks is described 
in \cite[section 10]{BoDe03}. The convergence of the diffusion
matrices in periodic environments is another illustration of these ideas.

The diffusion matrix is related to the matrix
$\A_0$  defined in \cite[section 5]{BoDe03} by homogenization.
In fact, $ \D_0=  2 \langle a(0)\rangle^{-1}\A_0$.
To see this, 
it suffices to note that for $\eta\in\R^d$, the function $f^\eta :\Omega\to\R^d$
defined by
$$f^\eta_i = \eta\cdot \chi_0(z_i) ,\quad 1\leq i\leq d$$
where $\{z_i ; 1\leq i\leq d\}$ is the canonical basis of $\R^d$,
is a solution of the auxiliary problem \cite[equation (15)]{BoDe03}.

By (\ref{lowerbdest}) and  (\ref{chinln}), $\pp$-a.s.
\begin{eqnarray*}
\Vert \dot  \chi_N\Vert_{\dot N}^2
   \leq  C\frac{1}{\dot \lambda_{N}} \tr\dot\E_N(\dot  \chi_N,\dot  \chi_N)
   \leq  C N^{d+2}.
\end{eqnarray*}
We have the following estimates  in terms of
the norm $\± \Vert u\Vert^2_{L^2(Q_N)}=\sum_{x\in Q_N}\vert u(x)\vert^2$ :
There is a constant $C<\infty$ 
such that for all $N\geq 1$,  $$ \Vert g\Vert_{L^2(Q_N)}^2  \leq   N^{d+2},\quad
\Vert\dot v_N\Vert_{L^2(Q_N)}^2 \leq C N^{d+2}\quad\hbox{and}\quad 
\tr\dot\E_N(\dot v_N,\dot v_N)\leq C N^d,$$
where $\±\dot v_N = \dot \chi_N +g.$

\bi
For  $\phii\in\C(\overline Q) $ or $\phii\in\C^\infty(\R^d) $, let $\± \hat \phii_N : Q_N \to\R$ be the function defined by $$\hat \phii_N(x) =\phii(x/N),\quad x\in Q_N.$$
By the diagonalization process, Riesz representation
theorem and 
  \cite[lemma 5]{BoDe03}, there is a subsequence $(N_k;k\geq 1)$ and
  a function $q\in\H^1(Q)^d$ such that 
for all $\phii\in(\C(\overline Q))^d $ and $f\in\{1, a(0,z_i), 1\leq i\leq d \}$
 as $N\to\infty$ along the subsequence $(N_k)$,
\begin{eqnarray}\label{one}
N^{-d-1}\sum_x\dot v_N(x)\cdot \hat\phii_N(x) f(x) & \to  & \langle f\rangle \int_Q q(s)\cdot \phii(s)ds
\end{eqnarray}
and for all $\± \phii\in(\C^\infty(\R^d))^d,$
\begin{eqnarray}\label{two}
N^{-d+1}\frac{1}{2}\dot\E_N( \dot  v_N , \hat\phii_N ) & \to  &  \int_Q 
(\nabla \phii)' \A_0 \nabla q ds 
\end{eqnarray}
where 
$\nabla \phii$ is the $d\times d$ matrix $(\nabla \phii)_{i,j}=\frac{\partial }{\partial s_j} \phii_i$,
$1\leq i,j\leq d$ and similarly for $\nabla q$.
By (\ref{one}),
\begin{eqnarray*}
N^{-d-1}\sum(\dot\chi_N + g)\cdot  \hat \phii_N 
 & = & N^{-d-1}\sum\dot\chi_N \cdot  \hat \phii_N  + N^{-d}\sum\frac{x}{N}\cdot   \hat \phii_N \\
&  \to &  0+\int_Q s\cdot \phii(s)ds.
\end{eqnarray*}

Therefore $\int_Q (q-g)ds = 0$
and $\tilde \chi_0\notation q-g$ is 1-periodic, that is $\tilde\chi_0\in\dot\H^1(Q)$
and $\int_Q \tilde \chi_0 ds=0$.

By (\ref{GG}), $\± \dot\E_N(\dot v_N ,\hat\phii_N) = 
(\dot H_N g,\hat\phii_N)_{\dot N} - (\dot H_N g,\hat\phii_N)_{\dot N} =0$.

Hence by (\ref{two}), $\± -\div \A_0\nabla (\tilde\chi_0 + g)=0$.
Then $\tilde \chi_0$ is the unique solution of the Poisson problem.
Therefore there is convergence in (\ref{one}) and (\ref{two})  for $N\to\infty$
and in particular, $\pp$-a.s.
\begin{eqnarray*}
\dot \D_N & = & \dot a_N(Q_N)^{-1}\dot\E_N(\dot v_N,\dot v_N)\\
 & = &  \dot a_N(Q_N)^{-1}\dot\E_N(\dot v_N, g)\\
   & \to & 2 \langle a\rangle^{-1}\ \int_Q (\nabla g)' \A_0\nabla q ds\\
   & = & 2 \langle a\rangle^{-1}\int_Q (\nabla g)' \A_0(\nabla g + \nabla \tilde\chi_0) ds\\
  & = & 2 \langle a\rangle^{-1}\A_0 =  \D_0.
\end{eqnarray*}
\qed

\section{Upper bounds on tail estimates}

\subsection{Martingale estimates and further notations}\label{martest}

The tail estimates
and the corresponding upper bounds on the variance
will be obtained by the method of bounded martingale differences
 developed by Kesten for first-passage percolation models  \cite{Kes93}.

When the conductances are assumed to be uniformly elliptic,
a stronger  version of Kesten's martingale 
inequality  can be used.
The same proof applies with some  simplifications.
In particular, \cite[Step (iii)]{Kes93} is not needed.
However the full generality of Kesten's martingale inequalities
 will be used when we consider  conductances
that are positive and bounded above but not 
necessarily uniformly elliptic. They are given in the second version.

 \bino\bf Martingale estimates I.\label{martest1}\rm\   
\smno\it
Let $\io (M_{n}; n\geq 0)$ be a martingale
on a probability space $(\Omega,\F,\pp)$.\hfill\break
Let $\Delta_{k}=M_{k}-M_{k-1}$, $k\geq 1$.
If there are positive random variables $U_k$
 (not necessarily $\F_k$-measurable)  such that
 for some  constant  $B_0<\infty$
\begin{equation}\label{martineq0}
  \EE(\Delta_k^2\vert \F_{k-1})\leq \EE (U_k \vert\F_{k-1})\   \hbox{for all} \    k\geq 1\quad
  \hbox{and}\quad  \sum_1^\infty U_k \leq B_0,
  \end{equation}
then 
$\io M_{n}\to M_{\infty}$ in $L^2$ and a.s., and  $\EE\vert M_\infty - M_0\vert^2 \leq 
B_{0}$.\sm
Moreover, if there is a constant $B_1<\infty$ such that for all $k\geq 1$,
\begin{equation}\label{martineq1}
\vert \Delta_k\vert \leq B_1,
\end{equation}
 then  for all $t>0$,
$$ \pp \left(\vert M_\infty - M_0\vert > t\right) \leq 4
\exp\left( -\frac{t}{ 4\sqrt B }\right)$$

where $\io B=\max\{B_0, e B_1^2\}$.
\rm

 \bino
 \bf Martingale estimates II.\label{martest2}\rm\   
\smno\it
Let $\io (M_{n}; n\geq 0)$ be a martingale
on a probability space $(\Omega,\F,\pp)$.\hfill\break
Let $\Delta_{k}=M_{k}-M_{k-1}$, $k\geq 1$.
If there is a constant $B_1<\infty$ such that for all $k\geq 1$,
\begin{equation}\label{martineqII1}
\vert \Delta_k\vert \leq B_1,
\end{equation}
if for some random variables $U_k\geq 0$
 (not necessarily $\F_k$-measurable)
 \begin{equation}\label{martineqII0}
  \EE(\Delta_k^2\vert \F_{k-1})\leq \EE (U_k \vert\F_{k-1})\   \hbox{for all} \    k\geq 1,
    \end{equation}
and if there are constants $0< C_1 ,C_2<\infty$ and  $s_0\geq e^2B_1^2$
such that,
\begin{equation}\label{martineqII2}
\pp(\sum_k U_k > s) \leq C_1 \exp(-C_2 s^2), 
\   \hbox{for all} \    s \geq s_0,
\end{equation}
then there are universal constants $c_1$ and $c_2$ 
which do not depend on $B_1$, $s_0$, $C_1$, $C_2$ nor on the distribution of $(M_k)$ and  $(U_k)$
such that   for all $s>0$,
$$ \pp\left(\vert M_\infty - M_0\vert > s\right) 
\leq c_1\left( 1+C_1+\frac{C_1}{s_0^2C_2}  \right)
\exp\left( -c_2\frac{s}{ s_0^{1/2} +  (s/(s_0C_2))^{1/3} }\right)$$
\rm

\bi
In the calculations to follow, we will use the following lighter notations.
In an environment $\omega$, the conductance of an edge $e$
with endpoints $x\sim y$ is denoted by $a(e,\omega)$
or $a(x,y,\omega)$. Similarly, 
for a function $v$ which is defined for $x$ and $y$, the endpoints of $e$,
the difference $v(x,\omega) - v(y,\omega)$ 
will be denoted, up to a sign, by $v(e,\omega)$.

Let $\±\L^d$ be the set of edges in $\Z^d$.
Whenever we assume that  the conductances are independent,
identically distributed
and bounded by $\kappa$, we will also assume that $\pp$
is a product measure on $\±\Omega = \lbrack 0, \kappa\rbrack^{\L^d}$
and $a(e,\cdot)$ are the coordinate functions.

Let $\{e_k ; k\geq 1 \}$ be a fixed ordering of $\±\L^d$
and  let $\F_k$,  $k\geq 1$, be the $\sigma$-algebra generated by
$\± \{ a(e_j,\cdot) ; 1\leq j\leq k\}$.
Then for an integrable random variable $h :\Omega \to \R$,
$$\EE(h  \vert   \F_k) = \EE_\sigma h(\lbrack \omega,\sigma\rbrack_k)$$
where $\lbrack \omega,\sigma\rbrack_k\in\Omega$ agrees with $\omega$ for the first $k$ coordinates and with $\sigma$ for all the other coordinates
and $\EE_\sigma$ denotes the integration with respect to $d\pp(\sigma)$.

\subsection{Approximations by a periodic environment}

In this section, we improve the tail estimates given in \cite{CaIo03}
by using  the regularity results from sections \ref{regres1} and \ref{regres2}
and the martingale estimates I  given in section \ref{martest1} above.
 Note that the calculations following the inequality \cite[(4.5)]{CaIo03}
 hold only for some laws $\pp$, like discrete 
 laws on a finite subset of $\R_+$.
 Denote the entries of $\dot D_N$ by $\dot\D_N^{ij}$, $1\leq i,j\leq d$.

\begin{thm}\label{thmappper}
Let $(a(e) ; e\in\L^d)$ be a sequence of i.i.d. uniformly elliptic conductances
with $\kappa$ as the ellipticity constant.
Then there is a constant $C$, $0<C<\infty$,
which depends only on the dimension and on the ellipticity constant $\kappa$
such that, 
 for all $t>0$ and  $N\geq 1$,
 \begin{equation}
    \label{eq:1}
     \max_{i,j}\pp(\vert \dot\D_N^{ij} - \EE \dot\D_{N}^{ij}\vert \geq t  N^{-\nu(d)} ) 
\leq 4\exp\left(- Ct\right),
 \end{equation}

\begin{equation}
    \label{eq:2}
    \max_{i,j}\vvar(\dot\D_N^{ij}) \leq 8C^{-2} N^{-2\nu(d)}
\end{equation}
where $\± 2\nu(d) = \max\left\{ \alpha, d-4+\alpha\right \}$ and
$\alpha>0$ is the regularity exponent which appears in (\ref{holdreg}).
 \end{thm}

\bi
Let $\{z_i ; 1\leq i\leq d\}$ be the canonical basis of $\R^d$.
For $1\leq i\leq d$, let $$\dot f_N = N^{-d}z_i'\dot\E_{N}(\dot 
v_{N},\dot v_{N}) z_i $$
where $$\dot v_{N}(x)= x+\dot\chi_{N}(x),\qquad x\in\overline Q_{N}.$$
\bi
Recall that by (\ref{convdifma}), $ \dot f_{N} \to \langle  a(0)\rangle z_i' \D_{0} z_i  $,
$\pp$ a.s. and in $L^1(\pp)$  as $N\to\infty$.

\bi
\begin{lem}\label{onedif}
Let $\omega$ and $\sigma$ be two environments such that
$a(e,\omega)= a(e,\sigma)$ for all edges $e$ except maybe for $e=e_k$.
Then for all $N\geq 1$,
$$\vert \dot f_N(\omega) -\dot f_N(\sigma)\vert  
\leq  \kappa  N^{-d}( \dot v^2_N(e_k,\omega) + \dot v^2_N(e_k,\sigma) )
$$
\end{lem}

\smno\it Proof of lemma  \ref{onedif}.\rm\  
By (\ref{poisseq}) and
by \ref{varpb1} of proposition \ref{poisson},
$\dot v_{N}$ is the solution of
a  variational problem.
 Then $\± \dot f_N(\omega) -\dot f_N(\sigma)$
\begin{eqnarray*}
 & =  &  N^{-d}\sum_e a(e,\omega)(z_i \cdot  \dot v_N(e, \omega) )^2  -  
  N^{-d}\sum_e a(e,\sigma) (z_i \cdot  v_N(e , \sigma) )^2 \\
 & \leq & N^{-d}\sum_e (\dot a_N(e,\omega) -\dot a_N(e,\sigma)) (z_i \cdot 
 \dot v_N(e,\sigma))^2\\
  & \leq  & \kappa N^{-d}   \dot v^2_{N}(e_k,\sigma).
\end{eqnarray*}
\qed

\bino\it Proof of theorem  \ref{thmappper}.\rm\  
Let $\io \Delta_{k} = \EE(\dot f_N\mid\F_k) - \EE(\dot f_N\mid \F_{k-1})$,
$k\geq 1$.
\smno
Let $M_0=0$ and $M_n =\sum_1^n \Delta_k$, $n\geq 1$.
Then $(M_n ; n\geq 0)$  is a martingale and we will see that
$$\dot f_N -\EE\dot f_N = \sum_1^\infty\Delta_{k}\qquad \hbox{a.s. and in 
$L^2$.}$$

We first check that  $(M_n ; n\geq 0)$ verifies conditions (\ref{martineq0})
and (\ref{martineq1}).

By (\ref{bdchin}), (\ref{bdchin24}) and the Hölder regularity (\ref{holdreg}),
there are constants $\beta$ and $C<\infty$ which depend only on $\kappa$ and $d$ 
such that $\pp$-a.s. and for all $N\geq 1$,
$$\sup_e   \dot v^2_{N}(e)  < C N^{\beta} 
$$
where $\± \beta = \min \{ d - \alpha , 4- \alpha \}$.
By lemma \ref{onedif}, we  see that
\begin{eqnarray*}
\vert\Delta_k\vert & = &\vert \EE(\dot f_N\mid \F_k) - \EE(\dot f_N \mid 
\F_{k-1})\vert\\
\\
 & = &  \vert\EE_\sigma(\dot f_N(\lbrack\omega,\sigma\rbrack_k )
 -\dot f_N(\lbrack\omega,\sigma\rbrack_{k-1}))\vert\\
 & & \\
 & \leq & \EE_\sigma \vert\dot f_N(\lbrack\omega,\sigma\rbrack_k ) 
-\dot f_N(\lbrack\omega,\sigma\rbrack_{k-1})\vert\\
 & & \\
 & \leq &  \kappa  N^{-d} \EE_\sigma(
 \dot v^2_N(e_k,\lbrack\omega,\sigma\rbrack_{k-1} ) +
\dot v^2_N(e_k,\lbrack\omega,\sigma\rbrack_{k}) )\\ 
\\
  & \leq &   C N^{-d+\beta}
\end{eqnarray*}
Hence, (\ref{martineq1}) holds with $B_1 =  C N^{-d+\beta}$.
Similarly,
\begin{eqnarray*}
\Delta_k^2 & \leq & \EE_\sigma(\vert\dot f_N(\lbrack\omega, \sigma\rbrack_k ) 
-\dot f_N(\lbrack\omega,\sigma\rbrack_{k-1})\vert^2)\\
 \\
 & \leq & 2\kappa^2 N^{-2d} \EE_\sigma(\dot v_N^4(e_k,\lbrack\omega,\sigma\rbrack_{k-1})  +
\dot v_N^4(e_k,\lbrack\omega,\sigma\rbrack_{k}))\\ 
  \\
 & \leq &C N^{-2d+\beta} \EE_\sigma(\dot v_N^2(e_k,\lbrack\omega,\sigma\rbrack_{k-1}) +
 \dot v_N^2(e_k,\lbrack\omega,\sigma\rbrack_{k})).
\end{eqnarray*}

Let $$ U_k(\omega) = 2C N^{-2d+\beta}\dot v_N^2(e_k,\omega).$$

We see that 
$\± \EE(\Delta_k^2 \mid \F_{k-1})\leq \EE(U_k\mid \F_{k-1})$ and that
 (\ref{martineq0}) holds with $B_0 =  2C N^{-d+\beta}$ since
\begin{eqnarray*}
\sum_k U_k  & =  &  2CN^{-2d+ \beta}  \sum_k\dot v_N^2(e_k,\omega)\\
 & \leq & C N^{-2d+ \beta} 
 \tr \dot\E(\dot v_N,\dot v_N)_{N,\omega}
 \leq   CN^{-2d+\beta+d} = C N^{-d+\beta}.
\end{eqnarray*}

Then the martingale estimates I hold with $B= \max\{ B_0 , e B_1^2\} = B_0$.
Hence,   for all $t>0$,
$$\pp(\vert\dot f_N -\EE\dot f_N\vert \geq t) \leq 
4\exp\left(-Ct N^{ (d-\beta)/2}\right)
$$
and for all $N\geq 1$, $\± \vvar(\dot f_N) \leq C N^{\beta - d }$.
\qed

\subsection{Effective conductance}

In this section, we obtain similar
tail estimates  for the effective conductances
of an increasing sequence
of cubes under mixed boundary conditions
and an upper bound on the variances.

It is simpler to work with these boundary conditions
because instead of the $L^\infty$ estimates
(\ref{bdchin})  and (\ref{bdchin24}),
we use the maximum principle :
if $u:\overline Q_N\to \R$ verifies $Hu = 0$ on $Q_N$ then
$\± \max_{\overline Q_N} u =\max_{\partial Q_N} u$.
Furthermore, since the maximum principle does not require
uniform ellipticity, it is possible to obtain good estimates
under weaker conditions on the conductances.
After the description of the model, we state two
theorems. The first one gives some  tail estimates when  
the conductances are uniformly elliptic  while the second one is  when they
are not.

Consider boundary conditions which can be interpreted as maintaining
a fixed potential difference   between two opposite faces of
$Q_{N}=\llbrack 1, N\rrbrack^d$
while the other faces are insulated. Denote the first coordinate of $x\in\Z^d$ by $x(1)$.
Let $\±   \V_{N}$ be the set of real-valued functions on $\overline Q_N$
such that
\begin{eqnarray*}
     u=0 \hbox{  on  }  
    \{x(1)=0\}\cap\partial Q_{N}, & & 
      u = N+1 \hbox{  on  }  
    \{x(1)=N+1\}\cap\partial Q_{N}\\
        \   \hbox{and} \quad u(x)=u(y) & &  \\
         \hbox{ for all } x\sim y , x\in \{x(1)\neq 
     0\}  &\cap & 
     \{  x(1)\neq N+1\}\cap\partial Q_{N}, y \in Q_{N}.
\end{eqnarray*}

 The Dirichlet form on $\H_N$, will be denoted by $\E_N$.
 For two functions $u,v:\overline Q_N\to\R$, it is defined by
  \begin{eqnarray}\label{dirform}
      \E_N(u,v) &  =  & 
\sum_{x,y} a(x,y)  (u(x)-u(y))(v(x)-v(y))
  \end{eqnarray}
 where  the sum is over all ordered pairs $\{x,y\}$
 such that $x\in Q_{N}$ and $y\in\overline Q_{N}$.

If for all edges $e$ of $\Z^d$, $a(e)>0$, then for all $N\geq 1$,
 there is a unique $\io v_{N} \in\V_N $ such that 
$Hv_{N}= 0$ in $Q_{N}$.
$v_{N}$ is also a solution of a variational problem :
it is the unique element of $\V_N$ such that
\begin{eqnarray*}
    \E_N(v_N,v_N) = \inf\big\{\   \E_N(u, u) ; u\in\V_N  \big \}.
\end{eqnarray*}
We will also write $\E_N(v_N,v_N)_\omega$ to indicate that
$\E_N $ and  $v_N$ are calculated with the conductances $a(e,\omega)$
 given by the environment $\omega$.

\sm
Let $\io f_N(\omega) = N^{-d}\E_{N}(v_{N}, v_{N})_\omega$.
It can be interpreted as 
the effective conductance between opposite faces of $Q_{N}$.

This model was considered by Wehr \cite{Weh97}. He
showed that for some laws, $\pp$, which include the 
exponential  and the one-sided normal distributions,
and assuming that $\EE (f_{N})$ is bounded below
by a positive constant, then $$ \liminf_{N} N^d \vvar( f_N )  > 0.$$

If the conductances are uniformly elliptic and stationary, 
then $\pp$ a.s. and in $L^1(\pp)$, as $N\to\infty$, 
$\io f_{N}$ converges to the effective conductance $f_0 = \nabla v_{0}\A_{0}\nabla v_{0}$
where $\A_{0}=\langle a(0)\rangle\D_0 $ is given in theorem \ref{statcov} and $v_{0}$
is the solution of a variational problem with mixed boundary conditions.
This was done in \cite[section 10]{BoDe03} by adapting
the homogenization methods of \cite[chapter 7]{JKO}.

\begin{thm}\label{effcond} 
Let $(a(e) ; e\in\L^d)$ be a sequence of i.i.d. uniformly elliptic conductances
on $\Z^d$, $d\geq 3$, with  ellipticity constant $\kappa\geq 1$. Let $C_0 = 32d\kappa^3$.
 
 Then for all $t\geq 0$ and  $N\geq 1$,
 $$\pp(\vert f_N -\EE f_N\vert \geq t N^{(2-d)/2}) 
 \leq 4\exp\left(- t/\sqrt{\kappa C_0}\right)
$$
and
$$\±   \vvar(f_N) \leq \kappa C_0 N^{2-d}. $$
\end{thm}

For conductances that are not necessarily uniformly elliptic, 
we have the following estimates.
Additional properties of non-uniformly elliptic reversible random
walks can be found in \cite{FoMa06}.

\begin{thm}\label{effcond2} 
Let $(a(e) ; e\in\L^d)$ be a sequence of i.i.d.  conductances
on $\Z^d$ such that for some constant
 $1\leq \kappa<\infty$,  $0< a(e)\leq \kappa$ for all edges $e$.
 Let $C_0 = 32d\kappa^3$.
\sm
 Then, for $d\geq 5$,
 $\± \vvar(f_N) \leq 128d\kappa^2 N^{4-d}$ 
 and  for all $t>0$ and  $N\geq 1$,
 \begin{equation}\label{firstest}
\pp(\vert f_N -\EE f_N\vert \geq t N^{(4-d)/2} )
\leq 4\exp\left(-  t/(8\kappa\sqrt {2d}) \right)
\end{equation}
\sm
If moreover, for some constants $D_0<\infty$ and $\gamma$, $ 0<\gamma < 2$, 
\begin{equation}\label{tailcond}
\pp( a^{-1}(e) \geq  s ) \leq D_0 s^{1-2/\gamma}, \quad\hbox{for all}\  s \geq 1,
\end{equation}
then, if $d\geq 4$ or if $d=3$ and $1/2\leq \gamma <1$, 
$$\vvar f_N \leq 16C_0(D_0+1)  N^{2-d+\gamma}$$
and for all $0< t < \left( \frac{C_0(D_0+1)^5}{8}\right)^{1/2} N^{(d-6+5\gamma)/2}$,
\begin{eqnarray}\label{tailpart}
\pp(\vert f_N - \EE f_N\vert > t) \leq 11c_1\exp\left(-\frac{c_2}{2\sqrt {2C_0(D_0+1)}}N^{(d-2-\gamma)/2}t \right).
\end{eqnarray}
where $c_1$ and $c_2$ are the constants that appear in the martingale estimates II.
In particular, they do not depend on $\kappa$, $d$ or $N$.

For $d=3$ and $0<\gamma \leq 1/2$,
$\± \vvar f_N \leq 16C_0(D_0+1)  N^{-1/2}$.
\end{thm}

\begin{lem}\label{lemconde}
Let $\omega$ and $\sigma$ be two environments such that
$a(e,\omega)= a(e,\sigma)$ for all edges $e$ except maybe for $e=e_k$.
Then for all $N\geq 1$,
$$\vert  f_N(\omega) -   f_N(\sigma)\vert  
\leq  \kappa  N^{-d}(  v^2_N(e_k,\omega) +  v^2_N(e_k,\sigma) )
$$
\end{lem}

\proof
Since $v_{N}\in\V_N$ is the solution of a variational problem,
\begin{eqnarray*}
f_N(\omega) -f_N(\sigma)
 & =  & N^{-d}( \E_N(v_N, v_N)_\omega -  \E_N(v_N, v_N)_\sigma) \\
  & \leq &   N^{-d}\sum_e (a(e,\omega) - a(e,\sigma)) v^2_N(e,\sigma)
 \leq \kappa N^{-d}v^2_N(e_k,\sigma)
\end{eqnarray*}
\qed

\bino\it
Proof of theorem \ref{effcond}.\rm\  
With the notations of \ref{martest},
let $\io \Delta_{k} = \EE(f_N\mid\F_k) - \EE(f_N\mid \F_{k-1})$,
 $M_0=0$ and $M_n =\sum_1^n \Delta_k$, $n\geq 1$.

To check that $(M_n ; n\geq 0)$ is a martingale that verifies conditions (\ref{martineq0})
and (\ref{martineq1}), 
 we have by lemma \ref{lemconde},
\begin{eqnarray*}
\vert\Delta_k\vert  & \leq & \EE_\sigma\left\vert f_N(\lbrack\omega,\sigma\rbrack_k ) 
- f_N(\lbrack\omega,\sigma\rbrack_{k-1})\right\vert\\
 & \leq &  \kappa N^{-d} \EE_\sigma v_N^2(e_k,\lbrack\omega,\sigma\rbrack_{k-1}) +
 \kappa N^{-d} \EE_\sigma v_N^2(e_k,\lbrack\omega,\sigma\rbrack_{k}) \\ 
  & \leq &  8\kappa N^{2-d}
\end{eqnarray*}
since by the maximum principle, $\± 0\leq v_N(x)\leq N+1$,
for all $x\in \overline Q_N$.

Then,
\begin{eqnarray*}
\Delta_k^2 & \leq & \EE_\sigma(\vert f_N(\lbrack\omega,\sigma\rbrack_k ) 
- f_N(\lbrack\omega,\sigma\rbrack_{k-1})\vert^2)\\
\\
 & \leq & 2\kappa^2 N^{-2d} \EE_\sigma v_N^4(e_k,\lbrack\omega,\sigma\rbrack_{k-1}) +
2\kappa^2 N^{-2d} \EE_\sigma v_N^4(e_k,\lbrack\omega,\sigma\rbrack_{k}) \\ 
 & & \quad\hbox{ and by the maximum principle,}\\
 & \leq & 8\kappa^2 N^{2-2d} \EE_\sigma v_N^2(e_k,\lbrack\omega,\sigma\rbrack_{k-1})  +
8\kappa^2 N^{2-2d} \EE_\sigma v_N^2(e_k,\lbrack\omega,\sigma\rbrack_{k}).\\
\end{eqnarray*}

For $k\geq 1$, let
\begin{equation}\label{sumuk}
 U_k(\omega) = 16\kappa^2 N^{2-2d} v_N^2(e_k,\omega).
\end{equation}

We see that $\± \EE(\Delta_k^2 \mid \F_{k-1})\leq \EE(U_k\mid \F_{k-1})$ and
by uniform ellipticity,
\begin{eqnarray*}
   \sum_k U_k    & = &  16\kappa^2 N^{2-2d}  \sum_k v_N^2(e_k,\omega)     \\
      &   &   \\
       & \leq  & 16\kappa^3 N^{2-2d} \E_N(v_N,v_N)_\omega \\
  & <  &  16\kappa^3(2d\kappa) N^{2-2d+d} = 32d\kappa^4 N^{2-d}.
\end{eqnarray*}

Hence  condition (\ref{martineq0}) holds with $B_0 = 32d\kappa^4 N^{2-d}$
and condition (\ref{martineq1}) holds with $B_1=4\kappa N^{2-d}$.
Therefore,
$\± f_N -\EE f_N = \sum_1^\infty\Delta_{k}$  a.s. and in  $L^2$
and for $d\geq 3$, by the martingale estimates I,
 we have that for all $N\geq 1$ and $t\geq 0$,
$$\pp(\vert f_N -\EE f_N\vert \geq t) \leq 
4\exp\left(-\frac{t}{4\sqrt B}   \right)
$$
where
$\± B = \max\{ B_0, eB_1^2\}
 =  32d\kappa^4 N^{2-d} $.
Accordingly, $\± \vvar(f_N) \leq 32d\kappa^4 N^{2-d}.$
\qed

\bino\it
Proof of theorem \ref{effcond2}.\rm\  
To obtain (\ref{firstest}), 
the preceding proof can be used up to (\ref{sumuk}) where
uniform ellipticity is first needed. 

Let $\± U_k(\omega) = 16\kappa^2 N^{2-2d} v_N^2(e_k,\omega)$, $k\geq 1$. 

Then simply bound $\± v_N^2(e_k,\omega)$ by $4N^2$ to obtain that
$$\sum_k U_k \leq 128d \kappa^2N^{4-d}.$$
Hence the martingale estimates I hold with
$\± B_0 = 128d \kappa^2 N^{4-d}, B_1=4\kappa N^{2-d}$
and $\± B=\max\{B_0, eB_1^2\}=B_0$.

These estimates can be improved  if we assume that (\ref{tailcond}) holds
 for some $0<\gamma<2$.
Starting from (\ref{sumuk}), we have that  for all $N\geq 1$,
\begin{eqnarray*}
\sum_k  v_N^2(e_k,\omega)  
 & \leq &  N^{\gamma}\sum_k a(e_k,\omega) v_N^2(e_k,\omega)
 + 4N^2 \sharp \{ k ; a(e_k,\omega)\leq N^{-\gamma}\}\\
   &\leq & 2d\kappa N^{d+\gamma}
  +  4N^{2}  \sharp \{ k ; a(e_k,\omega)\leq N^{-\gamma}\} . \\
\end{eqnarray*}
Therefore, $\±
\sum_k U_k
\leq  32d\kappa^3 N^{2-d+\gamma}
  + 64 \kappa^2 N^{4 - 2d}  \sharp \{ k ; a(e_k,\omega)\leq N^{-\gamma}\} .$
\smno
Let $C_0 = 32d\kappa^3$.
Then for  $t>1$,  $N\geq 1$ and $\gamma >0$,
\begin{eqnarray*}
\pp\left(\sum_k U_k > C_0t N^{2-d+\gamma}\right)
 & \leq & \pp\left ( \sharp \{ k ; a(e_k,\omega)\leq N^{-\gamma}\} >( t-1) N^{d-2+\gamma} \right) \\
  &\leq & 2\exp\left(-\frac{N^d}{4} ((t-1) N^{\gamma - 2} - p_N)^2 \right)
  \end{eqnarray*}
 by Bernstein's inequality with $\± p_N = \pp(a^{-1}(e)> N^\gamma)$.
 By (\ref{tailcond}),
 we have that for all $N\geq 1$,
 $\± p_N\leq D_0 N^{\gamma -2}$.
 Therefore, if $t> 2(1+D_0)$ then $(t-1-D_0)^2>t^2/4$ and
 \begin{eqnarray*}
\pp\left(\sum_k U_k > C_0t N^{2-d+\gamma}\right)
& \leq & 2\exp\left(-\frac{N^d}{4} (t-1-D_0)^2 N^{2\gamma - 4}  \right)\\
 & \leq  & 2\exp\left(-\frac{t^2}{16} N^{d-4+2\gamma}    \right)
  \end{eqnarray*}
 Equivalently, for all $\± s\geq 2C_0(D_0+1) N^{2-d+\gamma}$,
 $$\pp(\sum_k U_k > s)\leq  2 \exp\left( -\frac{s^2}{16C_0^2}N^{3d-8} \right).$$
 
 We see that the conditions for the martingale estimates II hold with
 $$B_1 = 4\kappa N^{2-d}, C_1=2,  C_2=\frac{1}{16C_0^2}N^{3d-8}\  \hbox{ and}\  s_0 =  2C_0 (D_0+1) N^{2-d+\gamma}
 $$
 since $C_0 \geq 8e^2\kappa^2$.

In particular,  we have  the tail estimates (\ref{tailpart}).

Moreover,
\begin{eqnarray*}
\vvar f_N & \leq &  \EE \sum_k U_k = \int_0^\infty \pp(\sum_k U_k> s )ds\\
 & \leq & s_0 + \int_{s_0}^\infty C_1 e^{-C_2 s^2}ds.
\end{eqnarray*}
Hence, if  $2-d+\gamma\geq 6-2d-\gamma$,
$$ \vvar f_N\leq  s_0 + \frac{C_1}{s_0 C_2}  \leq 16C_0(D_0+1) N^{2-d+\gamma} $$
while for $ d=3$
   and $ \gamma -1\leq -1/2$,
   $$\vvar f_N\leq  s_0 + \frac{C_1}{\sqrt C_2}     \leq  16C_0(D_0+1) N^{-1/2}.$$
\qed

\subsection{Spectral gap with Dirichlet boundary conditions}

In this last example,
we obtain tail estimates  for the spectral gap
of the random walk on  an increasing sequence
of cubes under Dirichlet boundary conditions.

Let   $\H_{N,0} =\{ u:\overline Q_N \to\R \  ;\  u=0\  \hbox{on}\  \partial Q_N \}$.
If $u,v \in \H_{N,0}$ then $\E_N(u,v) = (H u,v)_N$ where the Dirichlet form
 $\E_N$ is defined in (\ref{dirform}).
 
Let $\psi_{N}\in\H_{N,0} $ be the  solution of the variational problem :
\begin{eqnarray*}
    \E_N(\psi_N,\psi _N) = \inf\big\{\   \E_N(u, u) ; u\in\H_{N,0},\  \Vert u\Vert_{2,N} =1  \big \}.
\end{eqnarray*}

Then $\psi_N$ is unique (up to a sign) and  is an eigenfunction of $H$ acting on 
$\H_{N,0}$. Let  $\lambda_N>0$ be the corresponding eigenvalue.
It was shown in \cite{BoDe03}, by homogenization methods as in
Kesavan \cite{Kes79}, 
that  $N^2\lambda_N$ converges $\pp$-a.s. and in $L^1(\pp)$ as $N\to \infty$
to the Dirichlet eigenvalue of a second-order elliptic operator with constant coefficients.

The $L^\infty$ estimates of the eigenfunction
is provided by the De Georgi-Nash-Moser theory
(see  \cite[section 2.1]{Dav89} and 
\cite[chapter 11]{Chu97}) :
\smno\it 
Let $(a(e) ; e\in\L^d)$, $d\geq 3$,  be a (non-random) sequence of uniformly elliptic conductances.
Then there is a constant $C<\infty$ which depends
only on the dimension and on the ellipticity constant such that
if for $N\geq 1$, $\psi \in\H_{N,0}$ is a normalized eigenfunction of $H$, that is, for some
$\lambda >0$,  $\± H \psi  = \lambda \psi $
on  $Q_N$,   then 
for all  $x\in Q_N$,
\begin{equation}\label{regularityeig} 
\vert\psi (x)\vert\leq C \lambda^{d/4}.
\end{equation}
\rm\  
\sm
\begin{thm}\label{specgap}
Let $(a(e) ; e\in\L^d)$ be a sequence of i.i.d.
 uniformly elliptic conductances on $\Z^d$, $d\geq 3$.
 Let
 $$  f_N =  N^2\lambda_N =N^2\E_N(\psi_N , \psi_N).
 $$
Then there is a constant $C<\infty$ which depends only on $d$ and $\kappa$
such that for all $t>0$ and  $N\geq 1$,
$$\pp(\vert f_N -\EE f_N\vert \geq t N^{(2-d)/2}) \leq
 4\exp\left(-  t/C\right). 
$$
and
$$\± \vvar(f_N) \leq C N^{2-d}.$$
\end{thm}

\begin{lem}
Let $\omega$ and $\sigma$ be two environments such that
$a(e,\omega)= a(e,\sigma)$ for all edges $e$ except maybe at $e=e_k$,
an edge  with endpoints $x_k\sim y_k$ of $\Z^d$, where
$$a(e_k,\omega)\geq  a(e_k,\sigma).$$
Then for all $N\geq 1$,
\begin{eqnarray*}
\vert  f_N(\omega) -f_N(\sigma) \vert \leq  
  C ( N^2 \psi_N^2(e_k,\sigma) +
\psi_{N}^2(x_k, \omega) + \psi_{N}^2(y_k, \omega)).
\end{eqnarray*}
\end{lem}

\proof
  By the variational principle,
\begin{eqnarray*}
 \lambda_N(\omega) -\lambda_N(\sigma)
 & = & \E_N(\psi_{N}(\omega),\psi_{N}(\omega))_\omega - 
 \E_N(\psi_{N}(\sigma),\psi_{N}(\sigma) )_\sigma\\
  & \leq & \gamma \E_N(\psi_{N}(\sigma),\psi_{N}(\sigma) )_\omega
  -\E_N(\psi_{N}(\sigma),\psi_{N}(\sigma) )_\sigma\\
   & & \hbox{where}\  \gamma =\Vert\psi_N(\sigma)\Vert^{-2}_{2,N,\omega}\\
    & = & (\gamma -1)  \E_N(\psi_{N}(\sigma),\psi_{N}(\sigma) )_\omega\\
    & & + \  \E_N(\psi_{N}(\sigma),\psi_{N}(\sigma) )_\omega
  -\E_N(\psi_{N}(\sigma),\psi_{N}(\sigma) )_\sigma\\
   & \leq &  \E_N(\psi_{N}(\sigma),\psi_{N}(\sigma) )_\omega
  -\E_N(\psi_{N}(\sigma),\psi_{N}(\sigma) )_\sigma\\
  & & \hbox{since}\  \gamma \leq 1\\
   & \leq & (a(e_k,\omega) -a(e_k,\sigma))
   \psi_N^2(e_k,\sigma) \\
   & \leq & \kappa \psi_N^2(e_k,\sigma) .
\end{eqnarray*}

Similarly, by the variational principle,
\begin{eqnarray*}
\lambda_N(\sigma) -  \lambda_N(\omega) 
  & \leq &\overline \gamma \E_N(\psi_{N}(\omega),\psi_{N}(\omega) )_\sigma -
  \E_N(\psi_{N}(\omega),\psi_{N}(\omega) )_\omega \\
   & & \hbox{where}\  \overline \gamma =\Vert\psi_N(\omega)\Vert^{-2}_{2,N,\sigma}\\
    & = & (\overline \gamma -1)  \E_N(\psi_{N}( \omega),\psi_{N}(\omega) )_\sigma\\
    & & + \  \E_N(\psi_{N}(\omega),\psi_{N}(\omega) )_\sigma
  -\E_N(\psi_{N}(\omega),\psi_{N}(\omega) )_\omega\\
   & \leq &   (\overline \gamma -1)  \E_N(\psi_{N}( \omega),\psi_{N}(\omega) )_\sigma
  \end{eqnarray*}
since for all edges $ a(e,\omega)\geq a(e,\sigma)$.

Note that for all $u\in\H_N$,
$\± 0\leq \Vert u\Vert_{2,N,\omega}^2 -  \Vert u\Vert_{2,N,\sigma}^2
\leq a(e,\omega)(u^2(x_k) + u^2(y_k)).$
Then since $\± \psi_N(\omega)$ is a normalized eigenfunction
and the conductances are uniformly elliptic,
$$
0\leq  \overline\gamma -1  \leq C (\psi_{N}^2(x_k, \omega) + \psi_{N}^2(y_k, \omega)).
$$
Then,
$\± \lambda_N(\sigma) -  \lambda_N(\omega) 
\leq CN^{-2}(\psi_{N}^2(x_k, \omega) + \psi_{N}^2(y_k, \omega))$.
\qed

\bino
{\it Proof of theorem \ref{specgap}}.
As in the two preceding situations, we will verify  conditions (\ref{martineq0})
and (\ref{martineq1}) for
 $\io \Delta_{k} = \EE(f_N\mid\F_k) - \EE(f_N\mid \F_{k-1})$.
\begin{eqnarray*}
\vert\Delta_k\vert  & \leq & \EE_\sigma(\vert f_N(\lbrack\omega,\sigma\rbrack_k ) 
- f_N(\lbrack\omega,\sigma\rbrack_{k-1})\vert)\\
 & = & \EE_\sigma(\vert \cdots \vert\ ;\   a(e_k,\omega)\geq a(e_k,\sigma))
+  \EE_\sigma(\vert \cdots \vert\  ;\   a(e_k,\omega) < a(e_k,\sigma))\\
 & \leq & 
 C \EE_\sigma( N^2 \psi_N^2(e_k,\lbrack\omega,\sigma\rbrack_{k-1}) +
\psi_{N}^2 (x_k, \lbrack\omega,\sigma\rbrack_{k})
+ \psi_{N}^2(y_k, \lbrack\omega,\sigma\rbrack_{k}))\\
 & & 
 + C \EE_\sigma( N^2 \psi_N^2(e_k,\lbrack\omega,\sigma\rbrack_{k}) +
\psi_{N}^2 (x_k, \lbrack\omega,\sigma\rbrack_{k-1})\\
 &  & \hskip1cm +   \psi_{N}^2(y_k, \lbrack\omega,\sigma\rbrack_{k-1}))\\
  & \leq &  C N^{2-d}
\end{eqnarray*}
since by (\ref{regularityeig}), $\± \vert \psi_N(x)\vert \leq C N^{-d/2}$,
for all $x\in Q_N$.

Pursuing the above calculations, and using  (\ref{regularityeig}) again,
we find that
\begin{eqnarray*}
\Delta_k^2 & \leq & C N^{4-d}\EE_\sigma(  \psi_N^2(e_k,\lbrack\omega,\sigma\rbrack_{k-1}) +
\psi_N^2(e_k,\lbrack\omega,\sigma\rbrack_{k}))\\
 & & \quad + CN^{-d} \EE_\sigma(\psi_{N}^2(x_k, \lbrack\omega,\sigma\rbrack_{k}) 
+ \psi_{N}^2(y_k, \lbrack\omega,\sigma\rbrack_{k})\\
 & &\qquad  +  \psi_{N}^2(x_k, \lbrack\omega,\sigma\rbrack_{k-1}) 
+ \psi_{N}^2(y_k, \lbrack\omega,\sigma\rbrack_{k-1}))
\end{eqnarray*}

Let $\± U_k(\omega) = C N^{4-d} \psi_N^2(e_k,\omega )
+ C N^{-d} (\psi_{N}^2(x_k,\omega )  + \psi_{N}^2(y_k,\omega))$. 

We see that 
$\± \EE(\Delta_k^2 \mid \F_{k-1})\leq \EE(U_k\mid \F_{k-1})$
and
\begin{eqnarray*}
\sum_k U_k 
 & <  & C N^{4-d} \sum_k \psi_N^2(e_k,\omega )
+ C N^{-d} \sum_k (\psi_{N}^2(x_k,\omega )  + \psi_{N}^2(y_k,\omega))\\
& = & C N^{4-d} \E(\psi_N,\psi_N)_\omega + C N^{-d}\Vert \psi_N\Vert_{2,N}^2\\
  & \leq & C N^{4-d-2}+C N^{-d} \leq C N^{2-d}.
\end{eqnarray*}

Hence  (\ref{martineq0}) holds with $B_0= C N^{2-d}$ 
and   (\ref{martineq1}) holds with $B_1= C N^{2-d}$.
Since $\± B= \max\{B_0, eB_1^2\} \leq  C N^{2-d}$,
we obtain that  for all  $N\geq 1$ and $t>0$,
$$\pp(\vert f_N -\EE f_N\vert \geq t) \leq 
4\exp\left(-\frac{t}{C} N^{(d-2)/2}\right)
\quad\hbox{and}\quad\vvar(f_N) \leq C N^{2-d}.
$$
\qed

\section{A random walk in a random potential}

Let $(a(e) ; e\in\L^d)$, $d\geq 1$,  be a (non-random) sequence of positive conductances
and let $V :\Z^d\times\Omega \to \lbrack 0,\infty\lbrack$.
Add to the graph a vertex  $\ddag$  and an edge between $\ddag$ and each vertex
$x\in \Z^d$. In the environment $\omega$, its conductance  is given by
$$
a(x,\ddag,\omega)=a(x)V(x,\omega).
$$

The transition probabilites of  the random walk on $\Z^d\cup\{\ddag\}$
will be denoted  by $\breve p(x,y)$
to distinguish them from the transition probabilities of the reversible random walk on $\Z^d$ that
are given by $p(x,y)=a(x,y)/a(x)$. The former are defined by 
$ \breve p(\ddag,\ddag)=1$, $ \breve p(\ddag,x)=0$
and in the other cases, they are given by the conductances of the edges :
$$
\breve p(x,y) = \frac{p(x,y)}{1+V(x)},\quad
\hbox{and}\quad \breve p(x,\ddag)= \frac{V(x)}{1+V(x)}, \quad x,y\in\Z^d.
$$

The survival probability after each step is 
$\± (V(x)+1)^{-1}\notation e^{-\theta(x)}.$
Now assume that $(V(x) ; x\in\Z^d)$ is a sequence of i.i.d. nonnegative
random variables on $(\Omega,\F,\pp)$ that are not concentrated on zero. 

Let $T=\inf\{ k \geq 0; X_k = \ddag\} $.
Then the higher order transition probabilities are given by Feynman-Kac formula
\begin{eqnarray*}
\breve p(x,y,k) & = & E_x\left(\prod_{j=0}^{k-1}(V(X_j)+1)^{-1} ; X_k = y,  k< T\right)
\end{eqnarray*}
where $E_x$ is the expectation with respect to the reversible random walk on $\Z^d$.
Let $\breve P_x$ be the induced probability on the paths starting at $x$.

The Green function is defined by
$\breve G(x, y) =\sum_{k=0}^\infty \breve p(x,y,k)$.
A short calculation shows that,
since $V$ is not concentrated on zero, 
then for all $x,y\in\Z^d$, $d\geq 2$,
$\± \breve G(x,y)$ is a random variable with finite
moments of all order.

For a direction $x\in\Z^d$, $x\neq 0$,  and $ N\geq 1$,
let $$ f_{N}(x,\omega) = -N^{-1}\log\breve G(x,Nx,\omega).$$

We now prove a lower bound on the variance.
The analogue for first passage percolation is given in
\cite[(1.13)]{Kes93}.

\begin{prop}
If $(\theta(x), x\in\Z^d)$, $d\geq 2$,
is a  sequence of i.i.d. nonnegative random variables 
not concentrated on zero and such that
 $\EE(\theta(0)^2)<\infty$,
then for all $ x\in\Z^d$, $x\neq 0$,  and $N\geq 1$,
$$\vvar f_{N}(x) \geq N^{-2}\vvar \theta (0).$$
\end{prop}

\bino\bf Remark.\rm\   
In the particular case of constant conductances, that is $a(e)=1$ for all edges $e\in\L^d$,
Zerner  \cite{Zer98b} (see also \cite{Szn96} and \cite{WWM}) showed that
if $\EE \theta(0) <\infty$ then $f_{N}(x)$ converges $\pp$-a.s.
If, moreover,  $\EE(\theta(0)^2)<\infty$,
and if for $d=2$, there is $\underline \nu$ such that $\theta \geq \underline \nu >0$,
then by \cite[Theorem C]{Zer98b} there is a constant $C<\infty$ such that 
for all $ x\in\Z^d$, $x\neq 0$,  and $N\geq 1$,
$$ \vvar f_{N}(x) \leq C\frac{\vert x\vert}{N}.$$

\begin{proof}
For $x,y\in\Z^d$,
let 
$\tau_y =\inf\{ k\geq 0 ; X_k = y\}$ 
and $\tau_x^+ =\inf\{ k\geq 1 ; X_k = x\}$
with the convention that $\inf\emptyset=+\infty$.

Then, by conditioning on the time of last visit to $x$, we see that for $x\neq y$,
\begin{eqnarray*}
\breve G(x,y)  =  \breve P_{x}(\tau_y<T)\breve  G(y,y)
 = \breve G(x,x)\breve P_{x}(\tau_y <\tau_x^+ )\breve  G(y,y)
\end{eqnarray*}
and
$\± \log  \breve G(x,y) = \log  \breve G(x,x)+\log \breve P_{x}(\tau_y<\tau_x^+)
+ \log  \breve G(y,y).$
Since these are  decreasing functions of  $V$, 
by the FKG inequality (see \cite{Bar05} for a recent account), they are pairwise positively correlated. Then,
\begin{eqnarray*}
\vvar \log  \breve G(x,y)  & \geq & \vvar(\log \breve P_{x}(\tau_y <\tau_x^+)) \\
  & & +2 \ccov(\log  \breve G(x,x),\  \log \breve P_{x}(\tau_y<\tau_x^+))\\
 & & \hbox{     }  +  2 \ccov(\  \log \breve P_{x}(\tau_y<\tau_x^+), \log  \breve G(y,y) )\\
& \geq & \vvar(\log \breve P_{x}(\tau_y<\tau_x^+))\\
& = & \vvar\left(\log \sum_{z\sim x}   e^{-\theta(x)}p(x,z)\breve P_{z}(\tau_y<\tau_x^+)\right)\\
& = & \vvar(\theta(x)) +\vvar\left(\log \sum_{z\sim x} p(x,z) \breve P_{z}(\tau_y<\tau_x^+)\right)\\
 & \geq & \vvar(\theta(0)).
\end{eqnarray*}
\qed
\end{proof}

\bino
Laboratoire de Mathématiques UMR 6205, 
Université de Bretagne Occidentale,
6, avenue Le Gorgeu,
CS 93837,
29238 BREST Cedex 3
FRANCE
\hfill\break\noindent
E-mail: boivin@univ-brest.fr

\end{document}